\documentclass[a4paper,11pt]{article}
\usepackage{amsmath,amscd,amssymb}
\usepackage{graphicx,psfrag}

\def \vs {\vskip}
\def \hs {\hskip}
\def \noi {\noindent}

\def \oo {{\cal O}}

\def \a {{\alpha}}
\def \b {{\beta}}

\def \Ga {{\Gamma}}
\def \vp {\varphi}

\def \p {{\mathbb P}}
\def \Z {{\mathbb{Z}}}

\def \cC {{\mathfrak{C}}}

\def \fl {\rightarrow}
\def \ot {\otimes}

\def \pic {{\rm Pic}}
\def \pu {{\mathbb P}^1}
\def \opu {{\oo_\pu}}
\def \Mor #1#2{{\bf{Hom}}_{#1}(\pu,#2)}

\def \dm {{\textsc{Proof. }}}

\def \tha #1#2{\noi{\bf#1{\uppercase{\footnotesize{#2}}}}}

\newtheorem{theor}{\tha{T}{heorem}}[section]
\newenvironment{theo}{
  \begin{theor}\hs -0.2 cm {\bf .} ---  }
{  \end{theor}}

\newtheorem{propo}[theor]{\tha{P}{roposition}}
\newenvironment{prop}{
  \begin{propo}\hs -0.2 cm {\bf .} ---  }
{  \end{propo}}

\newtheorem{lemma}[theor]{\tha{L}{emma}}
\newenvironment{lemm}{
  \begin{lemma}\hs -0.2 cm {\bf .} ---  }
{  \end{lemma}}

\newtheorem{fait}[theor]{\tha{F}{act}}

\newtheorem{defini}[theor]{\tha{D}{efinition}}
\newenvironment{defi}{
  \begin{defini}\hs -0.2 cm {\bf .} ---  }
{  \end{defini}}

\newtheorem{corollaire}[theor]{\tha{C}{orollary}}
\newenvironment{coro}{
  \begin{corollaire}\hs -0.2 cm {\bf .} ---  }
{  \end{corollaire}}

\newtheorem{exemple}[theor]{\sc{Example}}
\newenvironment{exem}{
  \begin{exemple}\hs -0.2 cm {\bf .} ---  }
{  \end{exemple}}

\newtheorem{remarq}[theor]{\sc{Remark}}
\newenvironment{rema}{
  \begin{remarq}\hs -0.2 cm {\bf .} ---  }
{  \end{remarq}}

\newtheorem{preuve}{\sc{Proof}}

\def \comp {{\mathfrak{ne}}}
\def \l {{\lambda}}
\def \sca #1#2{\left\langle#1,#2\right\rangle}

\def \Xt {{\widetilde{X}}}
\def \ft {{\widetilde{f}}}

\def \at {{\widetilde{\a}}}

\def \fc {{\widehat{f}}}
\def \ac {{\widehat{\a}}}

\def \Del {\Delta}
\def \cln {{\cC_{L,n}(G/P)}}
\def \irin {{\mathfrak{ir}}}

\def \sp {{\underline{P}}}

\def \gq {{\underline{Q}}}
\def \sv {\underline{v}}
\def \su {\underline{u}}
\def \oq {\overline{Q}}
\def \sa {{\underline{a}}}
\def \sq {{\underline{q}}}
\def \zaq {Z(\sa,\sq)}

\def \p {{\mathbb P}}
\def \pu {{\mathbb P}^1}

\def \G {{\mathbb G}}

\def \fl {\rightarrow}
\def \ot {\otimes}
\def \vs {\vskip}
\def \hs {\hskip}
\def \l {{\cal L}}
\def \z {Z(\sp,\gq)}

\def \oo {{\cal O}}
\def \dm {{\textit{Proof --- }}}

\def \gg {{\mathfrak g}}
\def \gb {{\mathfrak b}}
\def \gp {{\mathfrak p}}
\def \gt {{\mathfrak t}}

\def \a {{\alpha}}
\def \b {{\beta}}
\def \gsp {{G/P}}
\def \zv {Z(\sp,\gq,\sv)}
\def \rv {R(\sp,\gq,\sv)}

\def \vp {\varphi}
\def \opu {{{\cal O}_{\pu}}}
\def \Si {{\Sigma}}
\def \pic {{\rm Pic}}

\def \noi {\noindent}

\def \U {\underline{{\cal U}}}
\def \V {\underline{{\cal V}}}

\def \Mor #1#2{{\bf{Hom}}_{#1}(\pu,#2)}

\def \Z {{\mathbb{Z}}}
\def \zuw {Z(\sp,\gq,\sw,\su)}

\def \th #1#2{\noi\textsc{#1 #2.} ---}

\def \C {{\mathbb C}}
\def \cC {{\mathfrak{C}}}
\def \comp {{\mathfrak{ne}}}
\def \oq {\overline{Q}}
\def \sa {{\underline{a}}}
\def \sq {{\underline{q}}}
\def \zaq {Z(\sa,\sq)}

\def \p {{\mathbb P}}
\def \pu {{\mathbb P}^1}

\def \G {{\mathbb G}}

\def \fl {\rightarrow}
\def \ot {\otimes}
\def \vs {\vskip}
\def \hs {\hskip}
\def \l {{\lambda}}

\def \z {Z(\sp,\gq)}

\def \oo {{\cal O}}
\def \dm {{\textsc{Proof. }}}

\def \gg {{\mathfrak g}}
\def \gb {{\mathfrak b}}
\def \gp {{\mathfrak p}}
\def \gt {{\mathfrak t}}

\def \a {{\alpha}}
\def \b {{\beta}}

\def \Ga {{\Gamma}}
\def \gsp {{G/P}}
\def \sca #1#2{\left\langle#1,#2\right\rangle}

\def \zv {Z(\sp,\gq,\sv)}
\def \rv {R(\sp,\gq,\sv)}

\def \vp {\varphi}
\def \opu {{{\cal O}_{\pu}}}
\def \Si {{\Sigma}}
\def \pic {{\rm Pic}}

\def \noi {\noindent}

\def \U {\underline{{\cal U}}}
\def \V {\underline{{\cal V}}}

\def \Mor #1#2{{\bf{Hom}}_{#1}(\pu,#2)}

\def \Z {{\mathbb{Z}}}
\def \zuw {Z(\sp,\gq,\sw,\su)}

\def \th #1#2{\noi\textsc{#1 #2.} ---}

\def \ne {N\!E}

\def \neff {{\rm Neff}}

\def \Xt {{\widetilde{X}}}
\def \ft {{\widetilde{f}}}

\def \at {{\widetilde{\a}}}

\def \fc {{\widehat{f}}}
\def \ac {{\widehat{\a}}}
\def \xaq {{X(\sa,\sq)}}
\def \F {{\mathbb{F}}}
\def \zaqk #1{Z(\sa-{#1},\sq-{#1})}
\def \sta {\widetilde{\sa}}
\def \stq {\widetilde{\sq}}

\def \Del {\Delta}

\begin{document}

~
\vs -0.5 cm

\centerline{\large{\uppercase{\bf{Rational curves on }}}}
\centerline{\large{\uppercase{\bf{homogeneous cones}}}}

\vs 1 cm

\centerline{\Large{Nicolas \textsc{Perrin}}}

\vs 1.5 cm

\centerline{\Large{\bf Introduction}} 

\vs 0.5 cm

In this text we study the scheme of morphisms from $\pu$ to any
homogeneous cone that is to say a cone $X$ over an homogeneous variety
$G/P$. Let us recall that we studied in \cite{PE} the scheme of
morphisms from $\pu$ to any homogeous variety. The main idea, in this
case is to restrict ourselves to the complementary of the vertex of
the cone, project on $G/P$ and apply the results of \cite{PE}.

\vs 0.2 cm

More precisely, let $G/P$ be an homogeneous variety and let $L$ be a
very ample divisor on $G/P$. We may embed $G/P$ in $\p(H^0L)$. If $V$ is
an $n$-dimensional vector space, it defines a linear subspace $\p(V)$
in $\p(H^0L\oplus V)$. Let us denote by $X=\cln$ the cone above
$G/P$ whose vertex is $\p(V)$. Now let $U$ be open subset of $X$
complementary to $Y=\p(V)$. We have a surjective morphism (see paragraph
\ref{preliminaires}):
$$s:\pic(U)^\vee\to A_1(X).$$
For any class $\a\in A_1(X)$, we can consider the following morphism:
$$i:\!\!\!\coprod_{s(\b)=\a}\Mor{\b}{U}\fl\Mor{\a}{X}$$ 
where $\Mor{\a}{X}$ is the scheme of morhisms $f:\pu\to X$ with
$f_*[\pu]=\a$ and $\Mor{\b}{U}$ is the scheme of morhisms $g:\pu\to U$
such that $[g]=\b$ where $[g]$ is the linear function $L\mapsto{\rm
  deg}(g^*L)$ on $\pic(U)$. As $Y=X\setminus U$ lies in codimension 2,
we expect the image of this morphism to be dense. For example we prove
in \cite{PEmin} that it is true for $X$ a minuscule Schubert variety.

\vs 0.2 cm

In our case the situation will be more complicated. Let us first
describe the "expected" components in the case where $i$ is
dominant. In this case we may apply the results of \cite{PE} to prove
that $\Mor{\b}{U}$ is irreducible as soon as it is non empty and the
images of these irreducible $\Mor{\b}{U}$ will give the irreducible
components of $\Mor{\a}{X}$. The expected components are thus indexed
by the subset $\comp(\a)$ of $\pic(U)^\vee$ given by elements
$\b$ such that $s(\b)=\a$ and $\Mor{\b}{U}$ is non empty. 

\vs 0.2 cm

This set can
be discribed in terms of roots: the ample divisor $L$ is a dominant
weight in the facet of the parabolic $P$. An element $\a\in A_1(X)$ is
completely determined by $\a\cdot L=d\in\Z$. Denote by
$\comp_B(\a)$ the set of all elements $\b$ in the cone 
generated by the positive roots such that $\sca{\b^\vee}{L}=d$. This
is a subset of $A_1(G/B)$. Then $\comp(\a)$ is its image in $A_1(G/P)$
(see paragraph \ref{preliminaires} for a more details). We prove the 

\begin{theo}
  Let $R$ be the
root lattice. 

{\rm (\i)} If $L(R)=\Z$, then the irreducible components of
the scheme of morphisms $\Mor{\a}{X}$ are indexed by $\comp(\a)$.

{\rm (\i\i)} If $L(R)\neq\Z$ (if we have $L> c_1(T_{G/P})$), then the
  irreducible 
  components of the scheme $\Mor{\a}{\mathfrak{C}(G/P)}$ are indexed by
  $\displaystyle{\coprod_{\a'\leq \a}\comp(\a')}$.
\end{theo}

We will see (paragraph \ref{preliminaires}) that $A_1(X)\simeq\Z$ so
that $\a'\leq\a$ in $A_1(X)$ means that the same inequality holds in
$\Z$. In the second case, we cannot deform a curve passing trough the
vertex of the cone so that the deformed curve does not pass trough it
any more. The integer $\a-\a'$ is then the multiplicity of the curve at
the vertex.

\begin{rema}
  {\rm (\i)} The condition $L(R)=\Z$ is exactly equivalent to the fact
  that there exists line on $G/P$ embedded with $L$. In other words
  there exists lines in the projectivized tangent cone to the
  singularity. 

We studied in \cite{PEmin} the same problem for minuscule Schubert
  varieties where the multiplicity in the singularity did not
  appear. If one consider more generaly quasi-minuscule Schubert
  varieties of non minuscule type (see \cite{GofG/P3} for a
  definition, the case of quasi-minuscule Schubert
  varieties of minuscule type should be very similar to the
  case of minuscule Scubert varieties) we
  recover this condition of the existence of lines in the
  projectivised tangent cone to the singularity. 

{\rm (\i\i)} If $P=B$ is a Borel subgroup and if we choose for $L$ the
Pl{\"u}cker embedding (or 
equivalently $L=\rho$ as a weight where $\rho$ is half the sum of the
positive roots) then $L(R)=\Z$ and the set $\comp(\a)$ is in bijection
with the set of irreducible integrable representations of level
exactly $\a\cdot L$ of the affine Lie algebra $\hat\gg$, see
paragraph \ref{preliminaires} for the general case).
\end{rema}

Here is an outline of the paper. In the first paragraph we define the
surjective map $s$ of the introduction and the set $\comp(\a)$ for an
homogeneous cone $X$. In the second paragraph we study the scheme of
morphisms from $\pu$ to the blowing-up $\Xt$ of the cone $X$ and prove
a smoothing result. In the last paragraph we prove our main result.

\vs 0.2 cm

The key point as indicated above is to study the surjectivity of the
map $i$ that is to say study the following problem: can any morphism
$f:\pu\to X$ be factorised in $U$ (modulo deformation). We do this by
lifting $f$ in $\ft$ on $\Xt$ and the problem becomes: does the lifted
curve $\ft$ of a general curve $f$ meet the exceptional divisor
$E$. If it is the case then we add a "line" $\Ga\subset E$ (this is
possible only when $L(R)=\Z$) with $\Ga\cdot E=-1$ and smooth the union
$\ft(\pu)\cup\Ga$. The intersection with $E$ is lowered by  one in the
operation. We conclude by induction on the number of intersection of
$\ft$ with $E$.

\vs 0.2 cm

We end with a discussion on the dimensions of the components, in
particular the variety $\Mor{\a}{X}$ is equidimensional if and only if
$L=\frac{1}{2}c_1(B/P)$ or $L=c_1(B/P)$.

\section{Preliminary}
\label{preliminaires}

In this paragraph we explain the results on cycles used in the
introduction. We describe the surjective morphism $s:\pic(U)^\vee\to
A_1(X)$ and define the set of classes $\comp(\a)$ for $\a\in A_1(X)$.

\vs 0.2 cm

Let $X$ be a scheme of dimension $n$. 
Denote by $Z_*(X)$ the group of 1-cycles on  $X$ and by
$Z^\equiv_*(X)$ and $Z^{r}_*(X)$ the subgroups of cycles trivial for
the numerical and rational equivalence. Let us denote by $N_*(X)$ 
and $A_*(X)$ the corresponding quotients. The Picard group is the
image in $A_{n-1}(X)$ of the subgroup of Cartier
divisors in $Z_{n-1}(X)$.

\begin{lemm}
  Let $X=\cln$ be a cone over a Schubert variety $G/P$ then

{\rm (\i)} $\pic(X)\simeq N^1(X)$,

{\rm (\i\i)}  $A_1(X)\simeq N_1(X)$. 

\noi
In particular we have $A_1(X)\simeq\pic(X)^\vee$ and they are
isomorphic to $\Z$. 
\end{lemm}

\dm
Consider with the decomposition $V\oplus H^0L$, the following group 
$$G'=\left(\begin{array}{cc}
GL(V)&{\rm Hom}(V,H^0L)\\
0&G
\end{array}\right)$$
acts on $X$ and the unipotent part $U(G')$ acts on $X$ with finitely
many orbits. In particular thanks to the results of \cite{FS..}

(\i) Thanks to the results of \cite{FS..} the groups $A_*(X)$ are free
generated by invariant subvarieties. The Picard group is contained in
$A_{n-1}(X)$ and is in particular free. Thanks to \cite{F}
Ex. 19.3.3. this implies that $\pic(X)\simeq N^1(X)$.

(\i\i) The results of \cite{FS..} also imply that $A_1(X)$ is
generated by the one-dimensional invariant subvarieties. The only such
subvarieties are the fibres of the cone so $A_1(X)\simeq\Z$ with a
fibre as generator. This fiber is clearly numerically free (for
example its degree is 1) so we get the result.

The duality comes from general duality between $N_1(X)$ and
$N^1(X)$.\hfill$\Box$

\vs 0.4 cm

Let $U$ be the smooth locus of $X$, it is also the dense orbit under
$G'$ in $X$. Let $Y$ be the complementary of $U$ in $X$, it is of
codimension at least 2 (at least when $\dim(G/P)>0$). This in
particular implies that $\pic(U)=A_{n-1}(U)\simeq A_{n-1}(X)$. We now
have the following inclusion: 
$$\pic(X)\subset A_{n-1}(X)\simeq\pic(U)$$
giving the surjection
$$s:\pic(U)^\vee\to A_1(X).$$
With these notations we make the following:

\begin{defi}
  Let $\a\in A_1(X)$. We define the set $\comp(\a)\subset
  A_{n-1}(X)^\vee$.

Let us make the
  identification $A_{n-1}(X)\simeq\pic(U)$. The elements of
  $\comp(\a)$ are the elements $\b\in\pic(U)^\vee$ such that
  $s(\b)=\a$ and there exists a complete curve $C\subset U$ with
  $[C]=\b$ as a linear form on $\pic(U)$ ($\b$ is effective).
\end{defi}

Let us describe $\comp(\a)$ explicitly: the smooth part $U$ is an
  affine bundle over $G/P$. In particular we have
  $\pic(U)\simeq\pic(G/P)$.

Let us Fix $T$ a maximal Torus in $G$, fix $B$ a Borel subgroup
containing $T$ and suppose that $B\subset P$. Let us denote by $\Del$ the set
of all roots, by $\Del^+$ (resp. $\Del^-$) the set of positive
(resp. negative) roots and by $S$ the set of simple roots associated to
the data $(G,T,B)$.

Denote by $\gg$, $\gt$ and $\gp$ the Lie algebras of $G$, $T$ and $P$
and define 
$$\a(\gp)=\Big\{\a\in S\ /\ \gg_\a\subset\gp\ {\rm and}\
\gg_{-\a}\not\subset\gp\Big\}.$$
Now set $\gt(\gp)^*$ as the subvector space of $\gt^*$ generated by the
roots in $\a(\gp)$, we have
$$\pic(G/P)\simeq\gt(\gp)\cap Q$$
where $\gt(\gp)$ is the dual of $\gt(\gp)^*$ in $\gt$ and $Q$ is the
weight lattice. The Picard group of $X$ in
$\pic(U)\simeq\pic(G/P)\simeq\gt(\gp)\cap Q$ is given by the
intersection of the line generated by $\l$ (the weight associated to
$L$) with the weight lattice $Q$. We have
$$\pic(U)^\vee\simeq\gt^*/\gt(\gp)^*\cap R$$
where $R$ is the root lattice.
Furthermore, an element $\b\in\pic(U)^\vee$ gives an effective element
if and only if it is in the image of the cone generated by positive
roots ie. in $\gt^*/\gt(\gp)^*\cap R^+$ (see \cite{PE}). Then we have
$$\comp(\a)=\left\{\b\in\gt^*/\gt(\gp)^*\cap R^+\ /\
  \sca{\b^\vee}{\l}=\a\cdot L\right\}$$
where the integer $\sca{\b^\vee}{\l}$ is well defined because
  $\l\in\gt(\gp)\cap Q$.

\begin{exem}
  Choose for $L$ (or for $\l$) the smallest ample sheave on
  $X$. This is possible: the picard group $\pic(U)=\gt(\gp)\cap Q$ is
  a direct sum of weight lattices of semi-simple Lie algebras
  $(\gg_i)_{i\in[1,r]}$. We just have to take
  $$\l=\sum_{i\in[1,r]}\rho_i$$ 
where $\rho_i$ is half the sum of positive roots in $\gg_i$.

Let us denote by $\irin_{\hat\gg_i}(\ell)$ the set of isomorphism classes of
irreductible integrable representations of level exactely $\ell$ of
the affine Lie algebra $\hat\gg_i$. Then we have
$$\comp(\a)=\prod_{\ell_1+\cdots+\ell_r=\a\cdot L}\irin_{\hat\gg_i}(\ell_i).$$
In particular if $P+B$ is a Borel subgroup of $G$ then $r=1$ and
$G_1=G$ and we recover the example of the introduction:
$$\comp(\a)=\irin_{\hat\gg}(\a\cdot L).$$
\end{exem}

\begin{rema}
  {\rm (\i)} The scheme $\Mor{\a}{X}$ is the scheme of
morphisms from $\pu$ to $X$ of class $\a$ (for more details see
\cite{GR} and \cite{MO}). 

In general, this will just mean that $\a\in A_1(X)$
  and that $f_*[\pu]=\a$ but sometimes (in particular in the
  introduction for the open part $U$) we consider $\a\in\pic(X)^\vee$
  and the class of a morphism $f:\pu\to X$ will be the linear form
  $\pic(X)\to\Z$ given by $L\mapsto{\rm deg}(f^*L)$.

In the case of a homogeneous cone $X$ the two notion coincide because
of the previous lemma.
In the case of the open part $U$ of a $X$, these scheme are connected
components of the scheme of morphisms with a fixed 1-cycle class
(which is always trivial).

{\rm (\i\i)} If $X$ is a variety, $\a\in A_1(X)$ and $F$ a vector bundle
on $X$ we will denote $\displaystyle{\a\cdot F=\int_\a c_1(F)}$ by
abuse of notation.
\end{rema}

\section{Resolution}

Recall that we denote by $X$ the cone $\cln$. Let $\Xt$ be the
blowing-up of $X$ in $\p(V)$. It is smooth and  
isomorphic to 
$$\p_{G/P}((V\ot\oo_{G/P})\oplus L).$$ 
Let us denote by $p$ the projection from $\Xt$ to $G/P$ and by
$\pi:\Xt\to X$ the blowing-up. The morphism $p$ has natural sections
given by points of $\p(V)$ or equivalently by surjective morphisms
$L\oplus(V\otimes\oo_{G/P})\to V\otimes\oo_{G/P}\to\oo_{G/P}.$

\subsection{Cycles on $\Xt$}

\begin{lemm}
  {\rm (\i)} Rational and numerical equivalences coincide on $\Xt$. In
  particular we have $A_1(\Xt)\simeq\pic(\Xt)^\vee\simeq
  A_{n-1}(\Xt)^\vee$.

{\rm (\i\i)} We have $\pic(\Xt)\simeq\pic(G/P)\oplus\Z$ with the
factor $\Z$ generated by the relative tangent sheaf $T_p$ of $p$.
\end{lemm}

\dm
(\i) Rational and numerical equivalence coincide on $G/P$. Moreover
the fibration in projetive spaces $\Xt\to G/P$ has sections so that
rational and numerical equivalences coincide on $\Xt$. This in
particular implies that $\pic(\Xt)=A_{n-1}(\Xt)=N^1(\Xt)$ and
$A_1(\Xt)=N_1(\Xt)$ and the duality follows.

(\i\i) The variety $\Xt$ is a $\p^n$-bundle over $G/P$ with sections
so we get that 
$$\pic(\Xt)\simeq\pic(G/P)\oplus\Z$$ 
with the factor $\Z$ generated by the relative tangent sheaf $T_p$ of
$p$.\hfill$\Box$

\vs 0.4 cm

Any element $\at\in A_1(\Xt)\simeq\pic(\Xt)^\vee$ is given by the
class $\b=p_*\at\in A_1(G/P)$ and the relative degree $d=\at\cdot
T_p$. We will use the notation $\ell=\b\cdot L=\at\cdot p^*L$.
Let us denote by $E$ the exceptional divisor on $\Xt$, it is a trivial
$\p^{n-1}$ bundle over $G/P$ given by the surjection
$L\oplus(V\otimes\oo_{G/P})\to V\otimes\oo_{G/P}$. Then we have: 
$$\at\cdot E=\frac{d-n\ell}{n+1},$$
it has to be an integer so that $d\equiv n\ell\ \ {\rm mod}\ n+1$. 

\vs 0.4 cm

Let us consider the following morphism still denoted $p$:
$$p:\Mor{\at}{\Xt}\fl\Mor{\b}{G/P}.$$

\begin{prop}
\label{irredenhaut}
  Thanks to the morphism $p$, the scheme $\Mor{\at}{\Xt}$ is an open
  subset of a projective bundle over $\Mor{\b}{G/P}$.
\end{prop}

\dm 
This generalises proposition 4 of \cite{PE} in the case where the
relative degree $\at\cdot T_p$ is negative. This is possible because
the vector bundle associated to the $\p^n$ fibration has a
decomposition $L\oplus(V\otimes\oo_{G/P})$.

Let $f:\pu\fl G/P$, we have to calculated the fiber of $p$ above
$f$. The fiber is given by
sections of the $\p^n$-bundle
$f^*(p):\p_{\pu}((V\otimes\opu)\oplus\opu(\ell))\fl\pu$ whose relative
degree is $d$. In other words the fiber is given by surjections
$(V\otimes\opu)\oplus\opu(\ell)\fl\opu(x)$ where $d=(n+1)x-\ell$
modulo scalar multiplication. The fiber is therefore isomorphic to an
open subset of $\p({\rm
  Hom}((V\otimes\oo_{\pu})\oplus\opu(\ell),\opu(x))$.

Let us remark that if $\Mor{\b}{G/P}$ is not empty then we
have $\ell\geq0$ and in this case $\Mor{\at}{\Xt}$ is not empty if and
only if $x\geq0$ when $n\geq2$ and if and only if $x=0$ or $x\geq\ell$
when $n=1$. In terms of $d$ this means that $d=-\ell$ or $d\geq n\ell$
if $n=1$ and $d\geq-\ell$ if $n\geq2$. In any cases, if
$\Mor{\at}{\Xt}$ is not empty then $x\geq0$.

There are two cases:
\begin{itemize}
\item If $x<\ell$ then any section is included in the exceptional
  divisor and the dimension of the fiber is:
  $$\frac{n}{n+1}(\ell+d)+n-1.$$
\item If $x\geq\ell$ then the fiber is of dimension $d+n$.\hfill$\Box$
\end{itemize}

Let $\at\in A_1(\Xt)$ such that $\Mor{\at}{\Xt}$ is not empty. This is
equivalent to the fact that $\b\in A_1(G/P)$ is positive (see
\cite{PE}, it is equivalent to the fact that $\Mor{\b}{G/P}$ is non
empty) and such that $d=-\ell$ or $d\geq n\ell$ if $n=1$, $d\geq -\ell$ if
$n\geq2$ (recall that $\ell=\b\cdot L$).

\begin{coro}
\label{dimension}
  The scheme $\Mor{\at}{\Xt}$ is irreducible of dimension
  \begin{itemize}
  \item $\displaystyle{\int_{\at}c_1(T_\Xt)+{\rm dim}(\Xt)}$ if $d\geq n\ell$
\item $\displaystyle{\int_{\at}c_1(T_\Xt)+{\rm
  dim}(\Xt)-\at\cdot E-1}$ if $d<n\ell$.
  \end{itemize}
\end{coro}

\dm 
We just use the preceding proposition and the fact proved in \cite{PE}
that the scheme $\Mor{\b}{G/P}$ is irreductible of dimension
$\int_{\b}c_1(T_{G/P})+{\rm dim}({G/P})$. Remark that in the last case
we have $n\ell>d$ so that the dimension of $\Mor{\at}{\Xt}$ is still
greater than the expected dimension
$\displaystyle{\int_{\at}c_1(T_\Xt)+{\rm dim}(\Xt)}$.\hfill$\Box$

\subsection{Smoothing curves on $\Xt$}

In this paragraph we will prove some results on curves on $\Xt$.

\begin{prop}
\label{lissage}
  Assume that $L(R)=\Z$. 

Let $\at\in A_1(\Xt)$, $\ft\in\Mor{\at}{\Xt}$ such that
  $\ft(\pu)\not\subset E$ and $\at\cdot E>0$. Assume that the
  image of $p\circ\ft:\pu\to G/P$ is not a line in the embedding given
  by $L$. 

Then there exists
  a deformation $\ft'$ of  $\ft$ and a curve $\Ga\subset\Xt$
  contracted by $\pi$ with
  $\Ga\cdot E=-1$ such that the curve $\ft'(\pu)\cup \Ga$ can be
  smoothed. The smoothed curve is the image of a morphism
  $\fc:\pu\to\Xt$.
\end{prop}

\dm
Let $(x,v)\in E\simeq G/P\times\p(V)$ be a point in the intersection
$\ft(\pu)\cap E$. 

\begin{lemm}
  There exists a deformation $\ft'$ of $\ft$ and a rational curve
  $\Ga$ in $\Xt$ such that $[\Ga]\cdot E=-1$, $[\Ga]\cdot L=1$ and
  meeting $\ft'(\pu)$ in exactly one point.
\end{lemm}

\dm 
Let us consider the lines in $G/P$ that is to say the rational curves
$\Ga'$ in $G/P$ such that $[\Ga']\cdot L=1$. Such curves exists
because we have $L(R)=\Z$. Let $\Ga'$ be such a line passing through
$p(x,v)=x\in G/P$ and let $\Ga$ be the section of $\Ga'$ in $E$ given
by the point $v\in\p(V)$. This curve is contracted by $\pi$ to the
point $v\in\p(V)$, its intersection with $E$ is given by $-[\Ga']\cdot
L=-1$.

As we assumed that $p\circ\ft(\pu)$ is not a line then $\Ga'$ meets
$p\circ\ft(\pu)$ in a finite number of points: $x$ and other points
$(x_i)$. The morphism $\ft$ is given by a section of the projective
bundle over $p\circ\ft$ that is to say by a surjection
$$s:(V\ot\opu)\oplus\opu(\ell)\fl\opu\left(\frac{d+\ell}{n+1}\right).$$
To deform $\ft$ we can deform this surjection $s$ in $s'$ such that at
$x$, we have $s'_x=s_x$ and at $x_i$ we have $s'_{x_i}\neq s_{x_i}$
for all $i$. This gives the required deformation.\hfill$\Box$  

\begin{lemm}
  The curve $\ft'(\pu)\cup\Ga$ can be smoothed. The smoothed curve is
  the image of a morphism $\fc:\pu\to\Xt$ of class $\ac$ with
$$\ac\cdot(p^*L+E)=\at\cdot(p^*L+E)\ \ {\rm and}\ \ \ac\cdot E<\at\cdot
E.$$ 
\end{lemm}

\dm 
If the smoothing exists then we have $\ac=\at+[\Ga]$ so $\ac\cdot
E=\at\cdot E-1$. Furthermore we have $p^*L+E=\pi^*L$ so that
$$\ac\cdot(p^*L+E)=\ac\cdot\pi^*L=\pi_*\ac\cdot
L=\pi_*\at\cdot L=\at\cdot\pi^*L =\at\cdot(p^*L+E).$$
This simply comes from the fact that $\pi_*[\Ga]=0$. 
Let us note that the curves $f'=\pi\circ\ft$ and $f''=\pi\circ\fc$ have
the same degrees but the curve $f''$ meets the vertex in one point less
than $f'$.

\vs 0.2 cm

To smooth $\ft'(\pu)\cup\Ga$ we use the following result proved in
\cite{HH} for $\p^3$ but valid for any smooth projective variety:

\begin{theo}
  Let $Z$ be a smooth projective variety and let $C$ be a nodal curve
  in $Z$. Assume that the cohomology group $H^1T_Z\vert_C$ is trivial
  then $C$ can be smoothed.
\end{theo}

We just have to prove that the cohomology group
$H^1(T_\Xt\vert_{\ft'(\pu)\cup\Ga})$ is trivial. We have the exact
sequences
$$0\fl T_{p}\fl T_\Xt\fl p^*T_{G/P}\fl 0\ \ {\rm and}\ \
0\fl\oo_{\ft'(\pu)}(-Q)\fl\oo_{\ft'(\pu)\cup
  \Ga}\fl\oo_\Ga\fl 0$$
where $Q$ is the intersection point of $\ft'(\pu)$ and $\Ga$. We just
have to prove the vanishing of the following cohomology groups:
$$H^1(p^*T_{G/P}\vert_{\Ga})\ \ ;\ \
H^1(p^*T_{G/P}\vert_{\ft'(\pu)}(-Q))\
\ ;\ \ H^1(T_{p}\vert_\Ga)\ \ {\rm and }\ \
H^1(T_{p}\vert_{\ft'(\pu)}(-Q)).$$
The first two groups are respectively equal to $H^1(T_{G/P}\vert_{\Ga'})$ and
$H^1(T_{G/P}\vert_{p(\ft'(\pu))}(-Q))$ where we denoted
$\Ga'=p(\Ga)$. They are trivial because $T_{G/P}$ is globally
generated and $\Ga'$ and $p(\ft'(\pu))$ are rational curves.

Let us denote by $\oo_{p}(1)$ the tautological quotient
of the projective bundle associated to $(V\ot\oo_X)\oplus L$, the
relative tangent sheaf is given by $T_{p}={\rm
  Coker}(\oo_\Xt\to((V^\vee\ot\oo_\Xt)\oplus
L^\vee)\ot\oo_{p}(1))$. In particular we have:
$$T_{p}\vert_\Ga={\rm Coker}(\opu\to(V^\vee\ot\opu)\oplus\opu(-1))\ \ \ \
\ \ {\rm and}\ $$
$$T_{p}\vert_{\ft'(\pu)}={\rm Coker}\left(\opu\fl
\left(V^\vee\ot\opu\left(\frac{d+\ell}{n+1}\right)\right)\oplus
\opu\left(\frac{d-n\ell}{n+1}\right)\right).$$ 
This proves that the group $H^1(T_{\pi}\vert_\Ga)$
vanishes. Furthermore, since $\ft$ exists we must have
$\frac{d+\ell}{n+1}\geq0$ (see proposition \ref{irredenhaut}) and
$\frac{d-n\ell}{n+1}=\at\cdot E>0$ so that
$H^1(T_{p}\vert_{\ft'(\pu)}(-Q))$ also vanishes.\hfill$\Box$

\section{Homogeneous cones}

Recall that we denote by $X$ the cone $\cln$. In this paragraph we
study the irreducible components of the scheme $\Mor{\a}{X}$ where
$\a\in A_1(X)$. Recall that $A_1(X)\simeq\Z$ and under this
identification $\a$ is just the degree of the corresponding curve.

\subsection{The case $L(R)=\Z$}

\begin{theo}
\label{sortir}
  Assume that $L(R)=\Z$, let $\a\in A_1(X)$ and
  $f\in\Mor{\b}{X}$. Then there exists a deformation $f'$ of $f$ such
  that $f'$ does not meet the vertex $\p(V)$ of the cone $X$.
\end{theo}

\dm 
Let us begin with the following:

\begin{lemm}
\label{facile}
  Let $f\in\Mor{\a}{X}$ such that $f$ factors through the vertex
  $\p(V)$ of the cone. Then there exists a deformation $f'$ of $f$ in
  $\Mor{\a}{X}$ such that $f'(\pu)$ does not factor through the
  vertex.
\end{lemm}

\dm 
Let $x\in G/P$ and consider the linear subspace generated by $x$ and
$\p(V)$. It is a $\p^{n+1}$ contained in $X$ and containing
$f(\pu)$. In this projective space we can deform the morphism $f$ so
that is does not factor through $\p(V)$ any more.\hfill$\Box$

\vs 0.4 cm

A general morphism $f\in\Mor{\a}{X}$ does not factor through the
vertex $\p(V)$ of the cone so it can be lifted in a morphism
$\ft:\pu\to\Xt$. Let $\at\in A_1(\Xt)$ the class of $\ft$, we have
$\pi_*\at=\a$. Because $f$ does not factor through the
vertex, the morphism $\ft$ does not factor through the
exceptional divisor $E$ so we have: $\at\cdot E\geq0$. If $\at\cdot
E=0$, then $\ft(\pu)$ does not meet $E$ thus $f$ does not meet the
vertex and we are done.
Let us assume that $\at\cdot E>0$. We proceed by induction on
$\at\cdot E$. Consider the morphism
$p\circ\ft:\pu\fl G/P$.

\begin{lemm}
  If the image of $p\circ\ft$ is a line in the projective embedding
  given by $L$ then there exists a deformation $f'\in\Mor{\a}{X}$
  of $f$ not meeting the vertex.
\end{lemm}

\dm 
Indeed, if the image of $p\circ\ft$ is a line then $f$ factors
through the linear subspace generated by the vertex and this line. It
is a $\p^{n+1}$ and the vertex is a linear subspace of codimension
2. There exists a deformation $f'$ of $f$ in this projective space not
meeting the vertex.\hfill$\Box$ 

\vs 0.4 cm

Let us now assume that the image of $p\circ\ft$ is not a line, we may
apply proposition \ref{lissage} so that there exists a deformation
$\ft'$ of $\ft$ and a curve $\Ga\subset\Xt$ contracted by $\pi$ with
$\Ga\cdot E=-1$ such that the curve $\ft'(\pu)\cup \Ga$ can be
smoothed. The smoothed curve is the image of a morphism
$\fc:\pu\to\Xt$ of class $\ac$. Let us consider $f'=\pi\circ\ft'$ and
$f''=\pi\circ\fc$. Then $f'$ is a deformation of $f$ and because $\Ga$
is contracted by $\pi$ the map $f''$ is a deformation of $f'$ and a
fortiori of $f$.

We have to prove the result on $f''$ whose lifting is $\fc$ of class
$\ac$. But we have $\ac\cdot E=\at\cdot E-1$ so the result is true by
induction.\hfill$\Box$

\begin{theo}
  Assume $L(R)=\Z$ and let $\a\in A_1(X)$ then the irreducible
  components of the scheme $\Mor{\a}{X}$ are indexed by
  $\comp(\a)$. For $\at\in\comp(\a)$ the dimension of the
  corresponding component is
$$\int_\at c_1(T_\Xt)+{\rm dim}(X).$$
\end{theo}

\dm 
Theorem \ref{sortir} proves that the set of morphisms $f:\pu\to X$
whose image does not meet the vertex $\p(V)$ is a dense open subset of
$\Mor{\a}{X}$. It is enough to study this open set. Any curve is this
open set comes from a unique lifting $\ft:\pu\to\Xt$ whose image does
not meet $E$. Let $\at\in A_1(\Xt)$ the class of $\ft$, since
$\at\cdot E=0$ we have $\at\in\pic(U)^\vee$ and in fact
$\at\in\comp(\a)$. The morphism
$$\pi_*:\coprod_{\at\in\comp(\a)}\Mor{\at}{\Xt}\to\Mor{\a}{X}$$
is thus dominant and birational (the inverse is given by lifting
morphisms). What is left to prove is that for each
$\at\in\comp(\a)$ the image of $\Mor{\at}{\Xt}$ (which is an
irreducible scheme) forms an irreducible component of
$\Mor{\a}{X}$. To prove this it is enough to prove that for any $\at$
and $\at'$ in $\comp(\a)$ the image of $\Mor{\at}{\Xt}$ is not
contained in the closure of $\Mor{\at'}{\Xt}$ in $\Mor{\a}{X}$. This
would be trivial if the scheme
$\coprod_{\at\in\comp(\a)}\Mor{\at}{\Xt}$ was equidimensional (it is
the case if $L=\frac{1}{2}c_1(G/P)$). In general, if it is the case
then there exists $f\in\Mor{\at}{\Xt}$ such that $f$ does not meet the
vertex and such that $f$ is the limit of a familly $f'_t$ of morphisms
in $\Mor{\at'}{\Xt}$. Because the condition of meeting the vertex is
closed me may assume that the elements $f'_t$ do not meet the
vertex. In particular projecting on $G/P$ gives a deformation from
$p(f'_t)$ to $p(f)$. This implies that $p_*\at=p_*\at'$ but as
$\at\cdot E=0=\at'\cdot E$ we have $\at=\at'$. The dimension comes
from corollary \ref{dimension}.\hfill$\Box$

\subsection{The case $L(R)\neq\Z$}

We begin with the following lemma on root systems:

\begin{lemm}
\label{equiv}
  Let $G$ be a semi-simple Lie group, $P\subset G$ a parabolic
  subgroup, $L$ a dominant weight in the facet defined by $P$ and $R$
  the lattice root, then we have the equivalence
$$L(R)\neq\Z\Longleftrightarrow L\geq c_1(G/P)$$
where $c_1(G/P)\in\pic(G/P)$ is considered as a weight and the order
is given by the positivity on simple roots.
\end{lemm}

\dm
Let us first describe $c_1(G/P)$ as a weight. Consider the set
$\a(\gp)$ of simple root and the lattice $\gt(\gp)\cap Q$ (which is
isomorphic to $\pic(G/P)$) defined in paragraph
\ref{preliminaires}. The lattice $\gt(\gp)\cap Q$ decomposes into a
direct sum of root lattices $R_i$. Let $\rho_i$ be half the sum of
possitive roots of the root system corresponding to $R_i$. Then we
have
$$c_1(G/P)=2\sum_i\rho_i.$$

If $L\geq c_1(G/P)$ then for any simple root $\a$ we have 
$$\sca{\a^\vee}{L}\geq\sca{\a^\vee}{c_1(G/P)}=
\sum_i\sca{\a^\vee}{\rho_i}=\left\{
  \begin{array}{cc}
0&{\rm if}\ \ \a\not\in\a(\gp)\\
2&{\rm if}\ \ \a\in\a(\gp)
  \end{array}\right.$$
and in particular $L(R)\subset2\Z$. 

Conversely, suppose that $L(R)\neq\Z$. Because $L$ is in the facet of
$P$ we have $\sca{\a^\vee}{L}=0$ for any simple root
$\a\not\in\a(\gp)$. If $\a$ is a simple root in $\a(\gp)$ then
$\sca{\a^\vee}{L}\geq2$ (otherwise $L(R)=\Z$). We see that for any
simple root $\sca{\a^\vee}{L}\geq\sca{\a^\vee}{c_1(G/P)}$ thus $L\geq
c_1(G/P)$.\hfill$\Box$

\begin{rema}
  Let $\at\in A_1(\Xt)$ such that $\at\cdot E\geq0$. Recall the
  notations $\b=\p_*\at$, $d=\at\cdot T_p$ is the relative degree and
  $\ell=\at\cdot p^*L=\b\cdot L$. Let $\a=\pi_*\at$ conidered as an
  integer. Then the dimension of $\Mor{\at}{\Xt}$ is given by 
$$\int_\at c_1(T_\Xt)+\dim(\Xt)=\int_\b
c_1(T_{G/P})+d+\dim(\Xt)=\int_\b c_1(T_{G/P})+(n+1)\at\cdot
E+n\ell+\dim(\Xt)$$
$$=\b\cdot (c_1(T_{G/P})-L)+(n+1)\at\cdot
(E+p^*L)+\dim(\Xt)$$
$$=\b\cdot (c_1(T_{G/P})-L)+(n+1)\a+\dim(\Xt).\ \ \ \ \ \ \ \ \ \ \ \
\ \! \ $$
So we have the formula
$$\dim(\Mor{\at}{\Xt})=\int_\at p^*(c_1(T_{G/P})-L)+(n+1)\a+\dim(\Xt).$$
\end{rema}

\begin{theo}
  Assume $L(R)\neq\Z$ and let $\a\in A_1(X)$. Then the irreducible
  components of $\Mor{\a}{X}$ are indexed by
  $\displaystyle{\coprod_{\a'\leq\a}\comp(\a')}$.
\end{theo}

\dm 
Thanks to lemma \ref{facile} (this lemma works without the hypothesis
$L(R)=\Z$) there exists a dense open subset of $\Mor{\a}{X}$ given by
morphisms $f$ that do not factor through the vertex of the cone. It is
enough to study this open set. In particular we know that the morphism
$$\pi_*:\coprod_{\at\in A_1(\Xt),\
  \pi_*\at=\a}\Mor{\at}{\Xt}\to\Mor{\a}{X}$$
is dominant. The classes $\at$ can even be choosen such that
  $\Mor{\at}{\Xt}$ is not empty. However the intersection $\at\cdot E$
  need not to be 0. In particular the classes $\at$ can be choosen in
$$A(\a)=\displaystyle{\coprod_{\a'\leq\a}\comp(\a')}$$
where $\a'=p_*\at\cdot L$ and $\a-\a'=\at\cdot E$. Indeed let $\at\in
A_1(\Xt)$ and set as usual $\b=p_*\at$. Then there exists a unique element
$\at'\in A_1(\Xt)$ such that $p_*\at'=\b$ and $\at'\cdot E=0$ (take
$n\b\cdot L$ for the relative degree). If $\at$ is such that
$\Mor{\at}{\Xt}$ is not empty then $\b$ is effective and because of
the value of the relative degree we have that $\Mor{\at'}{\Xt}$ is not
empty. In particular $\at'\in\comp(\a')$ for
$\a'=\pi_*\at'=p_*\at\cdot L$ and we have $\at\cdot
E=\pi_*\at-\pi_*\at'$. The element $\at$ is uniquely determined by
$\at'$ and $\at\cdot E$.

It is enough to prove that the images by $\pi_*$ of the irreducible
schemes $\Mor{\at}{\Xt}$ for $\at\in A(\a)$ are the irreducible
components. In other words we have to prove that for any $\at$ and
$\ac$ in $A(\a)$ the image of $\Mor{\at}{\Xt}$ is not contained in
the closure of the image of $\Mor{\ac}{\Xt}$ in $\Mor{\a}{X}$.

Let $\ft\in\Mor{\at}{\Xt}$ a generic point and $\fc_t$ a familly of
morphisms in $\Mor{\ac}{\Xt}$ such that $\pi\circ\fc_t$ converges to
$\pi\circ\ft$. In the compactification of $\Mor{\ac}{\Xt}$ by stable
maps (see for example \cite{FulPan}), the familly
$\fc_t$ has a limit say $\fc$ which is a morphism from a tree
$\cup_iD_i$ of rational curves to $\Xt$. Then we must have
$\pi\circ\fc=\pi\circ\ft$ as stable maps. In particular all but one of
the images by $\fc$ of the irreducible components of the tree are
contracted by $\pi$. To fix notation say that $D_i$ is contracted by
$\pi$ for $i\geq2$ and $\pi\circ\fc\vert_{D_1}=\pi\circ\ft$. Because
$\ft$ is generic, it is not contained in the exceptional divisor so
that the equality $\pi\circ\fc\vert_{D_1}=\pi\circ\ft$ implies that
$\fc\vert_{D_1}=\ft$. We see that $\fc_*[D_1]=\at$ so that
$$\ac=\fc_*[D_1]+\sum_{i\geq2}\fc_*[D_i]=
\at+\sum_{i\geq2}\fc_*[D_i].$$
In particular we have $\widehat{\b}=p_*\ac\geq p_*\at=\widetilde{\b}$
and because $L(R)\neq\Z$ we know thanks to lemma \ref{equiv} that
$L\geq c_1(G/P)$ and we get 
$$\widehat{\b}\cdot (c_1(T_{G/P})-L)\leq\widetilde{\b}\cdot
 (c_1(T_{G/P})-L).$$
As $\a=\pi_*\ac=\pi_*\at$ we see that 
$$\dim(\Mor{\ac}{\Xt})\leq\dim(\Mor{\at}{\Xt}).$$
But the morphism $\pi_*$ is generically injective on $\Mor{\ac}{\Xt}$
and $\Mor{\at}{\Xt}$ so that the scheme $\pi_*(\Mor{\at}{\Xt})$ cannot
be in the closure of $\pi_*(\Mor{\ac}{\Xt})$.\hfill$\Box$

\begin{rema}
  Let us end with a discussion on the dimension of the irreducible
  components of $\Mor{\a}{X}$ for $\a\in A_1(X)$. 

{\rm (\i)} In the first case: $L(R)=\Z$, these irreducible components
are indexed by elements $\at\in\comp(\a)$. For such an element we have
$\at\cdot E=0$ and the dimension of the component is given by
$$\dim(\Mor{\at}{\Xt})=\int_\at
p^*(c_1(T_{G/P})-L)+(n+1)\a+\dim(\Xt).$$
The "variable" part in this dimension is the first one and it is given
by
$$\b\cdot (c_1(T_{G/P})-L)$$
with $\b=p_*\at$ and we have $\a=\b\cdot L$ so that the "variable"
part is $\b\cdot c_1(T_{G/P})$. The element $\b$ ranges in the subset
of the positive cone in the root lattice $R$ (in the projection of $R$
in $\pic(G/P)$) given by the condition
$\b\cdot L=\a$. In particular if $L$ is not collinear to $c_1(G/P)$ the
dimensions of the irreducible components are not equal. In this case
the variety $\Mor{\a}{X}$ is equidimensional if anf only if
$L=\frac{1}{2}c_1(G/P)$. 

{\rm (\i\i)} In the second case: $L(R)\neq\Z$, these irreducible
components are indexed by elements
$\displaystyle{\at\in\coprod_{\a'\leq\a}\comp(\a')}$. For such an
element we have $\at\cdot E\geq0$ and the dimension of the component is given by
$$\dim(\Mor{\at}{\Xt})=\int_\at
p^*(c_1(T_{G/P})-L)+(n+1)\a+\dim(\Xt).$$
The "variable" part in this dimension is the first one and it is given
by
$$\b\cdot (c_1(T_{G/P})-L)$$
with $\b=p_*\at$. In this case we have $\b\cdot L=\a'\leq\a$. The
element $\b$ ranges in the subset of the positive cone in the root
lattice $R$ (in the projection of $R$
in $\pic(G/P)$) given by the condition $\b\cdot L\leq\a$. In particular if $L$
is not collinear to $c_1(G/P)$ the dimensions of the irreducible
components are not equal (look at the $\b$ such that $\b\cdot
L=\a$). Furthermore even if $L$
is not collinear to $c_1(G/P)$ the dimensions of the irreducible
components are not equal unless $L=c_1(G/P)$. In this case the variety
$\Mor{\a}{X}$ is equidimensional if anf only if $L=c_1(G/P)$.

\end{rema}

\begin{small}

\vs 0.2 cm

\noi
{\textsc{Institut de Math{\'e}matiques de Jussieu}}

\vs -0.1 cm

\noi
{\textsc{175 rue du Chevaleret}}

\vs -0.1 cm

\noi
{\textsc{75013 Paris,}} \hs 0.2 cm{\textsc{France.}}

\vs -0.1 cm

\noi
{email : \texttt{nperrin@math.jussieu.fr}}

\end{small}

\end{document}